\documentclass[a4j]{article} 

\usepackage{amssymb}

\title{Miura: Divisor Class Group Arithmetic}
\author{Joe Suzuki}
\date{}

\newtheorem{algo}{Algorithm}

\begin{document}

\maketitle


\begin{abstract}
The Package Miura contains functions that compute divisor class group arithmetic for nonsingular curves. The package reduces computation in a divisor class group to that in the ideal class group via the isomorphism. The underlying quotient ring  should be over the ideal given by a nonsingular curve in the form of Miura.
Although computing the multiplication of two integral ideals is not hard, we need to obtain an ideal such that the shortest Gr\"{o}bner basis component is the minimum
 in order to obtain the representative of the ideal class. Although the basic  procedure is due to Arita,  the source code has become much shorter using MaCaulay2.
The package is useful not just for computation itself but also for understanding the divisor class group arithmetic 
from the ideal point of view.

\end{abstract}

\section*{Introduction}

Let $K$ be a field, and $E:=\{(x,y)|y^2=x^3+ax+b\}$ with  $a,b\in K$ such that $4a^3+27b^2\not=0$.
It is known that if we define the arithmetic
$$P=(x_P,y_P)\in E \Longrightarrow -P=(x_P,-y_P)\in E$$
$$P=(x_P,y_P), Q=(x_Q,y_Q)\not=-P \Longrightarrow P+Q=(x_R,y_R)\in E\ s.t.$$
$x_R=\lambda^2-x_P-x_Q$, $y_R=\lambda(x_P-x_R)-y_P$, and 
$\displaystyle \lambda=\frac{y_Q-y_P}{x_Q-x_P}$ if $P\not=Q$ and $\displaystyle \lambda=\frac{3x_P^2+a}{2y_P}$ otherwise, then 
the set $E\cup \{P_\infty\}$ makes an addition group with zero element $P_\infty$.
Such a group arithmetic is currently being utilized in many applications such as elliptic curve cryptography and algebraic geometric coding theory.
The same notion is extended to non-singular curves with genus $g\geq 1$ in which
each group element is identified with a tuple of $g$ points.
Such a group is referred to as a divisor class group.

In order to realize such an arithmetic, the following should be considered:
\begin{enumerate}
\item how to express a curve as an Affine algebraic set like $E$ above for an elliptic curve. In particular, what if the curve is not
hyper-elliptic? A standard formula for the nonsingular curves could be used.
\item how to express each group element? In particular, if $g\geq 2$, more than one zero degree divisors exist for each class.
The representative element should be unique and computed efficiently.
\end{enumerate}
Those problems were solved in the late 1990's, and we will see the solutions in the following sections.

The package name ``Miura'' is due to the Japanese mathematician Shinji Miura \cite{miura} who gave the solution of the first problem.
The theory is not well known worldwide because his papers and Ph.D thesis were written only in Japanese.

The main feature of this package is that the source code is so simple. In fact, various solutions \cite{arita,basiri,hs} were proposed  to
the second problem. Although some of them run fast \cite{hs}, however, none of them are simple enough.

Macaulay2 \cite{m2} 
is a free computer algebra system developed by Daniel Grayson 
and Michael Stillman  for computation in commutative algebra and algebraic geometry.
This paper claims that using Macaulay2, the source code has become much simpler and 
the relation between the problems and solutions above has been much clear.
The procedure is essentially due to Arita \cite{arita} and source code has been improved by using Macaulay2.

The package, as of Dec 25, 2015, contains seven functions: pR, qR, inv, reduced, add, double, and multi (for the usage of those
functions will be explained in the later sections). For example, for the elliptic curves $y^2=x^3-3x$ over $\mathbb C$,
the addition of $(0,0)$ and $(1,2)$ is executed as follows.

\begin{verbatim}
i1 : R=pR(CC,{x,y},{2,3})
o1 = R
o1 : PolynomialRing
i2 : I=qR(R,y^2-x^3-3*x)
o2 = I
o2 : QuotientRing
i3 : J=ideal(x,y)
o3 = ideal (x, y)
o3 : Ideal of I
i4 : K=ideal(x-1,y-2)
o4 = ideal (x - 1, y - 2)
o4 : Ideal of I
i14 : add(J,K)
o14 = ideal (x - 3, y + 6)
o14 : Ideal of I
\end{verbatim}

\section*{Miura Theory}

Let $F$ be an algebraic function field of one variable over field $K$ \cite{st}.
If $P$ is the maximal ideal (place) of some variable ring $\cal O$ in  $F/K$, i.e.,
${\cal O}=\{z\in F|z^{-1}\not\in P\}$,
then the residue field $F_P:={\cal O}/P$ is defined for each $P\in {\mathbb P}_F$,
where ${\mathbb P}_F$ is the set of places in $F/K$. We define 
the the degree
$\deg P:=[F_P:K]$ and discrete variation $v_P: F\rightarrow {\mathbb Z}\cup \{\infty\}$
for each $P\in {\mathbb P}_F$, and
choose any  $P \in {\mathbb P}_F$ such that $\deg P=1$   (denote it by ${P}_\infty$),
and consider the vector space ${\cal L}(\infty P_\infty)$ over $K$:
$${\cal L}(\infty P_\infty):=\{f\in F|-v_P(f)\geq 0,P\not=P_\infty, P\in {\mathbb P}_F\}$$

Let $a_1,\cdots,a_t$ be generators of monoid ${\cal M}:=\{v_\infty(f)|f\in L\}$, and assume that 
they are positive integers and mutually prime.
Then, there exist $x_1,\cdots,x_t\in F$  such that $-v_\infty(x_1)=a_1,\cdots,-v_\infty(x_t)=a_t$.
Miura [2] considered that the kernel $ker \Theta$ of the surjective homomorphism
$\Theta: K[X_1,\cdots,X_t]\rightarrow K[x_1,\cdots,x_t]$
can be an affine algebraic set w.r.t. $F/K$, where $X_1,\cdots,X_t$ are indeterminates.


Let $\mathbb N$ be the nonnegative integers, and 
$\Psi:{\mathbb N}^t \rightarrow {\mathbb N}$ such that $\Psi(N):=\sum_{i=1}^tn_ia_i$ for $N=(n_1,\cdots,n_t)\in {\mathbb N}^t$.
We define the order $\prec$ among ${\mathbb N}^t$ by 
\begin{enumerate}
\item $\Psi(M)<\Psi(N) \Longrightarrow M\prec N$
\item $\Psi(M)=\Psi(N)$, and $m_1=n_1,\cdots,m_{j-1}=n_{j-1}, m_j>n_j$ for some $j=1,\cdots,t$ $\Longrightarrow M\prec N$
\end{enumerate}
for $M=(m_1,\cdots,m_t),N=(n_1,\cdots,n_t)\in {\mathbb N}^t$.
Moreover, we define another order $<$ among ${\mathbb N}^t$ by
$m_1\leq n_1,\cdots m_t\leq n_t$ and $m_j\not= n_j$ for some $j=1,\cdots,t$ $\Longrightarrow$
$M<N$ for $M=(m_1,\cdots,m_t),N=(n_1,\cdots,n_t)\in {\mathbb N}^t$.
Then, there exists  a minimum element in ${\mathbb N}^t$ (denote it by $b(m)$) w.r.t. $\prec$ such that $\Psi(N)=m$ for each $m\in {\cal M}$, 
and we define $B:=\{b(m)|m\in {\cal M}\}$.
Furthermore,  there exist minimal elements in ${\mathbb N}^t$
w.r.t. $<$ such that $N\not\in B$, and we define the set consisting of those elements as $V$.

For example, 
$B=\{(1,0),(0,1),(2,0),(1,1),\cdots\}$ and $V=\{(0,2)\}$ for $(a_1,a_2)=(2,3)$,
and 
$B=\{(1,0,0),(0,0,1),(0,1,0),(2,0,0),(1,0,1),
(1,1,0),(0,1,1),(3,0,0)\cdots\}$ and $V=\{(0,0,2),(0,2,0)\}$ for $(a_1,a_2,a_3)=(4,6,5)$.

Miura \cite{miura} proved the following statements:
\begin{enumerate}
\item each element $F_M$ in generator  $\{F_M|M\in V\}$ of ideal $ker \Theta$ is an element of 
\begin{equation}\label{eq1}
X^M+Span_K\{X^N|N\in B,\Psi(N)\leq \Psi(M)\}\backslash Span_K\{X^N|N\in B,\Psi(N)<\Psi(M)\}\ ,
\end{equation}
and satisfies
\begin{equation}\label{eq2}
Span_K\{X^N|N\in B\}\cap Ideal\{F_M|M\in V\}=\{0\}
\end{equation}
\item
If we are given mutually prime positive integers 
$a_1,\cdots,a_t$ and $\{F_M|M\in V\}$ that satisfies (1) and (2), then 
$Ideal\{F_M|M\in V\}$
is a prime ideal and the quotient field of $K[X_1,\cdots,X_t]/I$ is 
an algebraic function field of one variable over field $K$.
\end{enumerate}

For example, the affine algebraic sets 
for $(a_1,a_2)=(2,3)$ and $(a_1,a_2,a_3)=(4,6,5)$
are expressed by
\begin{equation}\label{eq3}
c_{0,0}+c_{1,0}x+c_{0,1}y+c_{2,0}x^2+c_{1,1}xy+x^3+y^2=0
\end{equation}
and 
\begin{equation}\label{eq4}
\left\{
\begin{array}{l}
c_{0,0,0}'+c_{1,0,0}'x+c_{0,0,1}'z+c_{0,1,0}'y+c_{2,0,0}'x^2+c_{1,0,1}'xz
+c_{1,1,0}'xy+z^2=0\\
c_{0,0,0}''+c_{1,0,0}''x+c_{0,0,1}''z+c_{0,1,0}''y+c_{2,0,0}''x^2+c_{1,0,1}''xz
+c_{1,1,0}''xy+c_{0,1,1}''yz++c_{3,0,0}''x^3+y^2=0
\end{array}
\right.
\end{equation}
using
and 
$c_{N}\in {\mathbb N}^2$ and 
$c_{N'}',c_{N''}''\in K$ with 
$N\in {\mathbb N}^2$ and $N',N''\in {\mathbb N}^3$, respectively.

Hereafter, we refer to the curves expressed by (1) and (2) as the Miura curves.

\section*{Divisor Class Group Arithmetic}
Hereafter, we assume that $K$ is a complete field.
Let ${\cal D}_F$, ${\cal P}_F$, and ${\cal C}_F:={\cal D}_F/{\cal P}_F$
 be the divisor group, principal divisor group, and divisor class group, respectively, of $F/K$, and 
define $\deg D:=\sum_{P\in {\mathbb P}_F}n_P \deg P$ for divisor $D=\sum_{P\in {\mathbb P}_F}n_PP \in {\cal D}_F$,
${\cal D}_F^0:=\{D\in {\cal D}_F|\deg D=0\}$, and ${\cal C}_F^0:={\cal D}^0_F/{\cal P}_F$.
We say that two divisors $D,D'\in {\cal D}_F$ are  equivalent, written $D\sim D'$, if $D=D'+(f)$ for some $f\in F\backslash \{0\}$, and
denote the class that contains  $D$ and its equivalent divisors w.r.t. $\sim$ by $[D]$.
We write $D=\sum_{P\not=P_\infty}n_PP\geq 0$ when $n_P\geq 0$ for $P\not=P_\infty$.

If $D\in {\cal D}^0$ is expressed as $E-nP_\infty$ with $E\geq 0$, then $D$ is said a semi-reduced divisor.
Moreover,
if $n$ is minimized in $D_1=E-nP_\infty$ with $E\geq 0$ and $D_1\sim D\in {\cal D}_F^0$, 
then $D_1$ is said the reduced divisor equivalent to the semi-reduced $D$. Then, it is known that
\begin{enumerate}
\item for each $C\in {\cal C}_F$, there exists a semi-reduced divisor $D$ such that $C=[D]$, and that
\item for each $D\in {\cal D}_F^0$, the reduced divisor $E-nP_\infty$ equivalent to $D$ is unique,
and $\deg E\leq g$, where $E\geq 0$.
\end{enumerate}
Thus, we can obtain reduced divisors using the following algorithm:
\begin{algo}\rm@
\begin{description}
\item[Input] Semi-reduced divisor $D=E-nP_\infty$ with $E\geq 0$
\item[Output] Reduced divisor $G\sim -D$
\end{description}
\begin{enumerate}
\item Find $f\in {\cal L}(\infty P_\infty)$ such that $(f)_0\geq E$ and the pole order $-v_{P_\infty}(f)$ is minimum, where 
$(f)_0:=\sum_{P: v(P)\geq>0}v_P(f)P$.
\item $G\leftarrow -D+(f)$
\end{enumerate}
\end{algo}
Since Algorithm 1 outputs a divisor equivalent to (-1) times the input divisor, 
a divisor equivalent to the input divisor can be obtained if Algorithm 1 is applied twice.
However, directly dealing with divisors is not generally efficient because of irreducible 
decomposition of polynomials.

Suppose that the curve expressed by (1) and (2) is nonsingular.
Then, the coordinate ring $K[x_1,\cdots,x_t]$ is Dedekind domain, so that
an isomorphism $\Phi$ from ${\cal C}_F$ to the ideal class  of $K[x_1,\cdots,x_t]$ is given as follows:
$$
[\sum_{P\not=P_{\infty}}n_PP-(\sum_{P\not=P_\infty}n_P)P_\infty]
\mapsto
[{\cal L}(\infty P_\infty - \sum_{P\not=P_{\infty}}n_PP)]\ ,
$$
where 
$${\cal L}(\infty P_\infty - \sum_{P\not=P_{\infty}}n_PP )
:=\{f\in F|-v_P(f)\geq n_P,P\not=P_\infty, P\in {\mathbb P}_F\}\ ,$$
and $[I]$ denotes the ideal class to which $I\subset K[x_1,\cdots,x_t]$ belongs. 
We say the ideals corresponding to reduced and semi-reduced divisors the reduced and semi-reduced ideals,
respectively, and each semi-reduced ideal $I$ is expressed by an integral ideal $I={\cal L}(\infty P_\infty - E)
\subset {\cal L}(\infty P_\infty)=K[x_1,\cdots,x_t]$ with $E\geq 0$.

Arita \cite{arita} proposed an algorithm to execute arithmetic on the Miura curves.
Note that the formula $-v_{P_\infty}(x^N)=\Psi(N)$ for  $N\in {\mathbb N}^t$
can be extended to $-v_{P_\infty}(f)=\max_N\Psi(N)$ for $0\not=f=\sum_N c_Nx^N\in {\cal L}(\infty P_\infty)$ with $c_N\in K$, $N\in {\mathbb N}^t$.
\begin{algo}\rm@
\begin{description}
\item[Input] Reduced ideals $J,K\in K[x_1,\cdots,x_t]$
\item[Output] The reduced ideal $N$ 
\end{description}
\begin{enumerate}
\item $L\leftarrow JK$
\item $f\leftarrow$ the minimum nonzero element in $L$ w.r.t. $\Psi$
\item $h\leftarrow$ the minimum nonzero element  w.r.t. $\Psi$ satisfying $(h)L\subset (f)$
\item $N\leftarrow (h/f)L$
\end{enumerate}
\end{algo}
In Algorithm 2, the minimum element in an ideal is computed by finding the reduced Gr\"{o}bner basis.
Although Algorithm 2 is simple enough, they proposed a detail implementation \cite{arita} of Algorithm 2.
For example, Maple requires many pages to implement the procedure \cite{basiri}.

For example, for $(a_1,a_2)=(2,3)$, if both  $(x,y)=(\alpha,\beta),(\alpha',\beta')\in K^2$ satisfy (\ref{eq3}), then
the minimum nonzero element in $L=JK$ for $J=(x-\alpha,y-\beta)$ and $K=(x-\alpha',y-\beta')$  w.r.t. $\Psi$
is the line connecting the two points:
$$f:\ (\beta'-\beta)(x-\alpha)-(\alpha'-\alpha)(y-\beta)=0\ .$$
If $f$ crosses with (\ref{eq3}) at $(\alpha'',\beta'')\in K^2$, then
then $(x-\alpha'',y-\beta'')$ is 
the quotient ideal $M:=(f)/L$. 
In a similar way,
the minimum nonzero element in $M=MP_\infty$ for $M=(x-\alpha'',y-\beta'')$ and $P_\infty=1$  w.r.t. $\Psi$
is the line connecting the two points:
$$h:\ x-\alpha''=0\ .$$
If $h$ crosses with (\ref{eq3}) at $(\alpha'',\beta''')\in K^2$, then
then $(x-\alpha'',y-\beta''')$ is 
the quotient ideal $N:=(h)/M$, where $\beta'''$ is known to be $-\beta''-a_{1,1}\alpha''-a_{0,1}$ \cite{suzuki}.


\section*{Implementation by MaCaulay2}

Our implementation using MaCaulay2 is pretty simple.
The first version (Dec 25, 2015) contains seven functions each of which consists of one sentence (eleven lines in total).

\begin{verbatim}
pR=(kk,v,w)->kk[v,MonomialOrder=>{Weights => w, Weights=>toList(1 ..#w) }];
qR=(R,p)->R/ideal p
inv=L-> quotient(ideal first first entries gens gb L, L)
reduced=L->inv inv L
add=(J,K)-> reduced (J*K)
double=J->add(J,J)
multi=(J,m)->(
    if m==0 then return ideal 1_I
    else if m\%2==0 then double multi(J,m//2)
    else add(double multi(J,(m-1)//2),J)
    )
\end{verbatim}

The function 
\verb+pR+ defines the polynomial ring $I$ and the weights $a_1,\cdots,a_t$ should be specified
as well as the constant field $K$ and variables $x_1,\cdots,x_t$.
For example, for curves over $\mathbb C$ and ${\mathbb F}_{5}$,  we specify them as in steps i1 and i16 below, respectively.

The function 
\verb+qR+ defines the quotient ring \verb+I=R/ideal p+ given the defined polynomial ring \verb+R+ and 
the ideal \verb+p+ of the algebraic set.
The argument \verb+p+ should be in the form (\ref{eq1}) 
and 
may be either an element in $R$ like (\ref{eq3}) or 
a list consisting of elements in $R$ like (\ref{eq4}).
 Although the condition (\ref{eq2})
is sometimes hard to check, it is known that the condition automatically is satisfied when 
the cardinality of $V$ is $t-1$ [5] (complete intersection).
For example, for elliptic and Miura curves expressed by $y^2=x^3+3x$ and
$\{y^2=x^3+1,z^2=xy+1\}$, we specify them as in steps i2 and i17 below, respectively.

The function \verb+inv+ computes the reduced ideal of $L^{-1}$ given a semi-reduced ideal $L$.
This implements Algorithm 1 in the context of ideal classes rather than divisor classes.
For example, for the \verb+L+ in step i5 which has been obtained via steps i3 and i4,
the reduced ideal of $L^{-1}$ is obtained via steps i6 through i12.
First of all,
applying gb to the ideal $L$ in step i6, 
we obtain an GroebnerBasis instance \verb+gb L+ in step o6.
Then, applying its method \verb+gens+ in i7, we obtain 
a matrix instance \verb+gens gb L+ in step o7.
Furthermore, applying \verb+entries+ in step i8 and \verb+first+ in steps i9 and i10,
we obtain the polynomial $f$: $y+5x=0$, the line  connecting $(0,0)$ and $(2,3)$, in o10.
Finally, the quotient ideal $(f)/L$ can be obtained in step o12.

The function 
\verb+reduced+ computes the reduced ideal of $L$ given semi-reduced ideals $L$
by repeating \verb+inv+ twice.

The function 
\verb+add+ computes the reduced ideal of $JK$ given semi-reduced ideals $J,K$ using the function \verb+reduced+.
For example, by applying \verb+inv+ once again to the ideal in step o12, we obtain the reduced ideal of $L$ in step o13.

The function \verb+double+ computes the reduced ideal of $L^2$
given a semi-reduced ideal $L$.

The function \verb+multi+ computes the reduced ideal of $J^m$ given semi-reduced ideals $J$
in a recursive way by applying \verb+add+ and  \verb+double+ at most $O(\log m)$ times:
$$J^m=
\left\{
\begin{array}{ll}
1&(m=0)\\
(J^\frac{m}{2})^2&(m: {\rm even})\\
(J^\frac{m-1}{2})^2J&({\rm otherwise})\\
\end{array}
\right.
$$
For example, we find  the reduced ideal of $K^6$ is a unit in step o15.

{
\begin{verbatim}
i1 : 
       -- Elliptic Curve
       R=pR(CC,{x,y},{2,3})
o1 = R
o1 : PolynomialRing
i2 : I=qR(R,y^2-x^3-3*x)
o2 = I
o2 : QuotientRing
i3 : J=ideal(x,y)
o3 = ideal (x, y)
o3 : Ideal of I
i4 : K=ideal(x-1,y-2)
o4 = ideal (x - 1, y - 2)
o4 : Ideal of I
i5 : L=J*K
             2                          3
o5 = ideal (x  - x, x*y - 2x, x*y - y, x  - 2y + 3x)
o5 : Ideal of I
i6 : gb L
o6 = GroebnerBasis[status: done; S-pairs encountered up to degree 2]
o6 : GroebnerBasis
i7 : gens gb L
o7 = | y-2x x2-x |
             1       2
o7 : Matrix I  <--- I
i8 : entries gens gb L
                2
o8 = {{y - 2x, x  - x}}
o8 : List
i9 : first entries gens gb L
               2
o9 = {y - 2x, x  - x}
o9 : List
i10 : first first entries gens gb L
o10 = y - 2x
o10 : I
i11 : ideal first first entries gens gb L
o11 = ideal(y - 2x)
o11 : Ideal of I
i12 : quotient(ideal first first entries gens gb L, L)
o12 = ideal (x - 3, y - 6)
o12 : Ideal of I
i13 : reduced(L)
o13 = ideal (x + 3, y + 2)
o13 : Ideal of I
i14 : add(J,K)
o14 = ideal (x - 3, y + 6)
o14 : Ideal of I
i15 : multi(K,6)
o15 = ideal 1
o15 : Ideal of I
\end{verbatim}
}

The same procedure is applied to the Miura curves than contain elliptic, hyperelliptic, $C_{ab}$ curves as special cases.
If the genus is $g$, each ideal class contains $g$ zeros 
$$L_1=(x_1-\alpha_{1,1},\cdots,x_t-\alpha_{1,t}),\cdots,L_g=(x_1-\alpha_{g,1},\cdots,x_t-\alpha_{g,t})$$
and the reduced ideal is obtained by reducing $L_1 \cdots L_g$.
In the following case, since the genus is four, the reduced ideal with zeros $J,K,L,M$ in steps i18,i19,i20,i21
is obtained by multiplying them  and reducing it as in step i22. 

{
\begin{verbatim}
i16 : 
      -- Miura Curves
      R=pR(GF 5,{x,y,z},{4,6,5});
i17 : I=qR(R,{y^2-x^3-1,z^2-x*y-1});
i18 : J=ideal(x-2,y-2,z);
o18 : Ideal of I
i19 : K=ideal(x-4,y,z-1);
o19 : Ideal of I
i20 : L=ideal(x,y-1,z-4);
o20 : Ideal of I
i21 : M=ideal(x,y-4,z-1);
o21 : Ideal of I
i22 : A=reduced(J*K*L*M)
              2
o22 = ideal (x  + y + z + 2x, x*z - 2y - 2z + 2x, x*y - y - z - x, y*z - 2y - 2z - x + 1)
o22 : Ideal of I
i23 : multi(A,654)
o23 = ideal 1
o23 : Ideal of I
i24 : multi(A,327)
o24 = ideal (x + 1, y)
o24 : Ideal of I
i25 : add(A,inv A)
o25 = ideal 1
o25 : Ideal of I

\end{verbatim}
}

\end{document}